\documentclass[leqno,12pt]{article}

\usepackage{amsfonts,amsmath,amssymb}
\usepackage{amsthm}
\usepackage{enumerate}
\usepackage{graphics}
\usepackage{graphicx}
\usepackage{float}%
\usepackage{enumerate}%
\usepackage{amsfonts,amsmath,amssymb}
\usepackage{amsthm}
\usepackage{enumerate}
\usepackage{graphics}
\usepackage{graphicx}
\theoremstyle{definition}
  \newtheorem{theorem}{Theorem}[section] %
  \newtheorem{corollary}[theorem]{Corollary}
  \newtheorem{defin}[theorem]{Definition}%
  \newtheorem{example}[theorem]{Example}%
  \newtheorem{problem}[theorem]{Problem}
  \newtheorem{basicproblem}[theorem]{Basic Problem}
  \newtheorem{remark}[theorem]{Remark}%
  \newtheorem{observation}[theorem]{Observation}
\newcommand{\invHom}[3]{\operatorname{Hom}_{#1}({#2},{#3})}
\numberwithin{equation}{section}
\begin{document}

\title
{
Topics on global analysis of manifolds
 and representation theory of reductive groups
}
\date{}                                           
% Activate to display a given date or no date
%
\author{
Toshiyuki Kobayashi
\\
Graduate School of Mathematical Sciences, 
The University of Tokyo, 
\\
and Kavli IPMU (WPI)
}
\maketitle %
{\abstract{
Geometric symmetry induces symmetries of function spaces, 
 and the latter yields a clue to global analysis 
 via representation theory.  
In this note we summarize recent developments
 on the general theory 
 about how geometric conditions affect representation theoretic properties
 on function spaces, 
 with focus on multiplicities and spectrum.  
}}

{\bf{Mathematic Subject Classification (2020):}}\enspace
Primary 22E46;
Secondary 43A85, 22F30
%%%%%%%%%%%%%%%%%%%%%%%%%%%%%%%%%%%%%%%%%%%%%%
\section{\lq\lq{Grip strength}\rq\rq\ of representations on global analysis --- geometric criterion for finiteness of multiplicities}
\label{sec:grip}
%%%%%%%%%%%%%%%%%%%%%%%%%%%%%%%%%%%%%%%%%%%%%%
To which extent, 
 does representation theory provide a useful information for global analysis
 on manifolds?

As a guiding principle, 
 we begin with the following perspective 
(\cite{xkronsetsu2019}).

\begin{basicproblem}
[{\lq\lq{Grip strength}\rq\rq\ of representations}]
\label{prob:1.1}
{\rm{
Support that a Lie group $G$ acts on $X$.  
Can the space of functions on $X$ be 
 \lq\lq{sufficiently controlled}\rq\rq\ by the representation 
theory of $G$?
}}
\end{basicproblem}

The vague words, \lq\lq{sufficiently controlled}\rq\rq,
 or conversely,
 \lq\lq{uncontrollable}\rq\rq,
need to be formulated as mathematics.  
Let us observe 
 what may happen in the general setting
 of infinite-dimensional representations
 of Lie groups $G$.

\begin{observation}
For an infinite-dimensional $G$-module $V$, 
 there may exist infinitely many different irreducible subrepresentations.  
Also, the multiplicity of each irreducible representation 
can range from finite to infinite.

When confronting such a general situation, one focuses on the principle:
\smallskip
\par
\quad 
$\bullet$\enspace
even though there are infinitely many (sometimes uncountably many)
 different irreducible representations,
the group action {\bf{can}} distinguish the different parts;
\smallskip
\par\quad
$\bullet$\enspace
the group action {\bf{cannot}} distinguish the parts where the same irreducible 
representations occur with multiplicities.  
\end{observation}

This observation suggests us to think of the multiplicity
 of irreducible representations
 as an obstruction of \lq\lq{grip strength of a group}\rq\rq.
For each irreducible representation $\Pi$ of a group $G$,
 we define the multiplicity of $\Pi$ in the regular representation $C^{\infty}(X)$ by
\begin{equation}
\label{eqn:multfn}
  \operatorname{dim}_{\mathbb{C}} \operatorname{Hom}_G
  (\Pi,C^{\infty}(X)) \in {\mathbb{N}} \cup \{\infty\}.  
\end{equation}
The case 
$
   \operatorname{dim}_{\mathbb{C}} \operatorname{Hom}_G
  (\Pi,C^{\infty}(X))=1
$
 (multiplicity-one) provides
 a strong  \lq\lq{grip strength}\rq\rq\
 of representation theory on global analysis, 
 which may be illuminated
 by the following example:
\begin{example}
\label{ex:PDE}
Let $X$ be a manifold, 
 $D_1, \cdots, D_k$ differential operators on $X$, 
 and $G$ the group of diffeomorphisms $T$ of $X$
 such that $T \circ D_j=D_j \circ T$
 for all $j=1, \cdots, k$.  
Then the space of solutions $f$
 to the differential equations on $X$:
\begin{equation}
\label{eqn:Joint}
  D_j f = \lambda_j f 
  \qquad
  \text{for $1 \le j \le k$}
\end{equation}
forms a $G$-module
 (possibly, zero)
 for any $\lambda_1, \cdots, \lambda_k \in {\mathbb{C}}$.  
The group $G$ becomes a Lie group
 if $\{D_1, \cdots, D_k\}$ contains the Laplacian
 when $X$ is a Riemannian manifold 
 (or more generally a pseudo-Riemannian manifold).  
Assume now that the multiplicity of an irreducible representation $\Pi$
 of $G$ in $C^{\infty}(X)$ is one.  
Then any function belonging to the image
 of a $G$-homomorphism from $\Pi$ to $C^{\infty}(X)$
 satisfies a system of differential equations \eqref{eqn:Joint}
 for some $\lambda_1, \cdots, \lambda_k \in {\mathbb{C}}$
 by Schur's lemma.  
\end{example}

We formalise Basic Problem \ref{prob:1.1} as follows.

\begin{problem}
[Grip strength of representations on global analysis]
\label{prob:6.2}
{\rm{Let $X$ be a manifold on which a Lie group $G$ acts.
Consider the regular representation of $G$ on $C^{\infty}(X)$ 
 by 
\[
C^{\infty}(X) \ni f(x) \mapsto f(g^{-1} \cdot x) \in C^{\infty}(X)
\quad
\text{for $g \in G$.}
\]

\smallskip
\par\noindent
(1)\enspace
Find a necessary and sufficient on the pair $(G,X)$
 for which the multiplicity \eqref{eqn:multfn} of every irreducible representation $\Pi$ of $G$ in the 
regular representation $C^{\infty}(X)$ is {\bf{finite}}.
\smallskip
\par\noindent
(2)\enspace
Determine a condition on the pair $(G,X)$ for which
the multiplicity is {\bf{uniformly bounded}} with respect to all irreducible representations $\Pi$.
}}
\end{problem}

A solution to Problem \ref{prob:6.2} will single out
 a nice setting of $(G,H)$
 in which we could expect a detailed study 
 of global analysis
 on the homogeneous manifold $X=G/H$
 by using representation theory of $G$.  
The multiplicity may depend on the irreducible representations $\Pi$ in (1),
 and thus we may think that the group $G$ has \lq\lq{stronger grip power}\rq\rq\ in (2)
 than in (1).  
We may also consider a {\bf{multiplicity-free case}}:
\par\noindent
(3)\enspace
Determine a condition
 on the pair $(G,X)$
 for which the multiplicity \eqref{eqn:multfn}
 is either 0 or 1
 for any irreducible representation $\Pi$ of $G$.

Clearly, 
 (3) is stronger than (2), 
 however, 
 we do not discuss (3) here.

Problem \ref{prob:6.2} is settled 
 in Kobayashi--Oshima \cite{xktoshima}
 for homogeneous spaces $X$ of reductive Lie groups $G$.  
To state the necessary and sufficient condition,
we recall some notions from the theory of transformation groups.
The following terminology was introduced in \cite{Ksuron}.

\begin{defin}
[Real sphericity]
\label{def:6.3}
Suppose that a reductive Lie group $G$ acts continuously 
on a connected real manifold $X$.
We say $X$ is a {\bf{real spherical}}
 if a minimal parabolic subgroup of $G$ has an open orbit in $X$.  
\end{defin}

As is seen in Example \ref{ex:sph} below, 
 the classical notion of 
 spherical varieties is a special case
 of real sphericity
 because a minimal parabolic subgroup
 of a complex reductive group $G_{\mathbb{C}}$ is nothing but
 a Borel subgroup.

\begin{example}
[Spherical variety]
\label{ex:sph}
Suppose that $X_{\mathbb{C}}$ is a connected complex manifold
 and that a complex reductive Lie group $G_{\mathbb{C}}$ acts 
biholomorphically on $X_{\mathbb{C}}$.  
Then $X_{\mathbb{C}}$ is called a {\bf{spherical variety}} of $G_{\mathbb{C}}$, 
 if a Borel subgroup of $G_{\mathbb{C}}$ has 
an open orbit in $X_{\mathbb{C}}$.  
Spherical varieties have been extensively studied
 in algebraic geometry, 
 geometric representation theory, 
 and number theory.  
\end{example}

Here are some further examples.  
\begin{example}
\label{ex:spherical}
Let $X$ be a homogeneous space of a reductive Lie group $G$ and 
$X_{\mathbb{C}}$ its complexification.
\par\noindent
{\rm{(1)}}\enspace
The following basic implications hold
(Aomoto, Wolf, and Kobayashi--Oshima).  

\begin{align*}
&\text{$X$ is a symmetric space}
\\
&\qquad\Downarrow \,\, \text{Aomoto, Wolf}
\\
&\text{$X_{\mathbb{C}}$ is a spherical variety}
\\
&\qquad\Downarrow \,\,\text{Kobayashi--Oshima \cite[Prop.~4.3]{xktoshima}}
\\
&\text{$X$ is a real spherical variety}
\\
&\qquad\Uparrow \,\,\text{obvious}
\\
&\text{$G$ is compact.}
\end{align*}

\par\noindent
{\rm{(2)}}\enspace
When $X$ admits a $G$-invariant Riemannian structure, 
 the following are equivalent
 (see Vinberg \cite{Vinberg2001}, Wolf \cite{Wolf2007}):
\begin{align*}
&\text{$X_{\mathbb{C}}$ is spherical}
\\
\Longleftrightarrow
&\text{$X$ is weakly symmetric in the sense of Selberg} 
\\
\Longleftrightarrow
&\text{$X$ is a commutative space.}
\end{align*}

\par\noindent
{\rm{(3)}}\enspace
The classification of irreducible symmetric spaces was accomplished by 
Berger \cite{ber} at the level of Lie algebras.  
\par\noindent
{\rm{(4)}}\enspace
The classification theory of spherical varieties $X_{\mathbb{C}}$
 has been developed by Kr{\"a}mer, 
Brion,
 Mikityuk, 
 and Yakimova.  
\par\noindent
{\rm{(5)}}\enspace
The triple space 
$(G \times G \times G)/\operatorname{diag} G$
is not a symmetric space.  
It is real spherical 
 if and only if $G$ is locally a direct product 
 of compact Lie groups
 and $S O(n,1)$, 
 see \cite{Ksuron}.  
This geometric result implies a finiteness criterion 
 of multiplicities
 for the tensor product
 of two infinite-dimensional irreducible representations
 (\cite{Ksuron}, \cite[Cor.~4.2]{xkProg2014}).  
The triple space is considered as a special case
 of the homogeneous space
 $(\widetilde G \times G) / \operatorname{diag} G$
 for a pair of groups
 $\widetilde G \supset G$.  
More generally, 
 the classification of real spherical manifolds
 $(\widetilde G \times G) / \operatorname{diag} G$
 was accomplished in \cite{xKMt} 
 when $(\widetilde G, G)$ are irreducible symmetric pairs
 in connection to the branching problem
 for $\widetilde G \downarrow G$, 
 see \cite{xkProg2014}.  
\par\noindent
{\rm{(6)}}\enspace
Let $N$ be a maximal unipotent subgroup
 of a real reductive Lie group $G$.  
Then $G/N$ is real spherical, 
 as is easily seen from the Bruhat decomposition.  
Moreover, 
 the following equivalence holds:
\[
\text{$G_{\mathbb{C}}/N_{\mathbb{C}}$ is spherical
$\Longleftrightarrow$ $G$ is quasi split.  }
\]
This is related to the fact
 that the theory of Whittaker models 
 ({\it{e.g.}} Kostant--Lynch, H.~Matumoto)
 yields
 stronger consequences
 when $G$ is assumed to be quasi split, 
 see Remark \ref{rem:Wh} below.  
\end{example}

We denote by $\operatorname{Irr}(G)$ the set
 of equivalence classes
 of irreducible admissible smooth representations
 of $G$.  
We do not assume unitarity for here.  
The solutions of Problem \ref{prob:6.2}, which is a 
reformalisation of Basic Problem \ref{prob:1.1}, 
are given by the following two theorems.

\begin{theorem}
[Criterion for finiteness of multiplicity
{\cite{xktoshima}}]
\label{thm:HP}
Let $G$ be a reductive Lie group and
$H$ a reductive algebraic subgroup of $G$,
and set $X=G/H$.
Then the following two conditions on the pair $(G,H)$
are equivalent.
\par\indent
{\rm{(i)}}\enspace
{\rm{(representation theory)}}\enspace
$\dim_{\mathbb{C}} \operatorname{Hom}_{G}(\Pi, C^{\infty}(X))< \infty$
 $({}^{\forall} \Pi \in {\operatorname{Irr}}(G))$.  
\par\indent
{\rm{(ii)}}\enspace
{\rm{(geometry)}}\enspace
$X$ is a real spherical variety.
\end{theorem}

In \cite{xktoshima}, 
 the proof of the implication 
 (ii) $\Rightarrow$ (i)
 uses (hyperfunction-valued) boundary maps
 for a system of partial differential equations
 with regular singularities, 
 whereas that of the implication (i) $\Rightarrow$ (ii)
 is based on a generalization
of the Poisson transform.  
This proof gives
 not only the equivalence of (i) and (ii) in Theorem \ref{thm:HP} 
but also some estimates of the multiplicity from above and below.  
In turn, 
 these estimates bring us
 to the following geometric criterion of the uniform boundedness
 of multiplicity.

\begin{theorem}
[Criterion for uniform boudedness of multiplicity {\cite{xktoshima}}]
\label{thm:HB}
Let $G$ be a reductive Lie group and
$H$ a reductive algebraic subgroup of $G$,
and set $X=G/H$.
Then the following three conditions on the pair $(G,H)$
are equivalent.
\par\indent
{\rm{(i)}}\enspace
{\rm{(representation theory)}}\enspace
There exists a constant $C$ such that 
\[
\dim_{\mathbb{C}} \operatorname{Hom}_{G}(\Pi, C^{\infty}(X)) \le C
\quad
({}^{\forall} \Pi \in {\operatorname{Irr}}(G)).
\]
\par\indent
{\rm{(ii)}}\enspace
{\rm{(complex geometry)}}\enspace
The complexification $X_{\mathbb{C}}$ of $X$ is 
a spherical variety of $G_{\mathbb{C}}$. 
\par\indent
{\rm{(iii)}}\enspace
{\rm{(ring theory)}}\enspace
The ring of $G$-invariant differential operators on $X$ is commutative.
\end{theorem}

\begin{remark}
{\rm{
The equivalence (ii) $\Leftrightarrow$ (iii) in Theorem \ref{thm:HB}
 is classical, 
 see {\it{e.g., }} \cite{Vinberg2001}, 
 and the main part here is to characterize
 the representation theoretic property (i)
 by means of conditions in other disciplines.  
}}
\end{remark}

\begin{remark}
{\rm{
In general, 
 the constant $C$ in (i) cannot be taken 
to be 1
 when $H$ is noncompact.  
}}
\end{remark}

\begin{remark}
{\rm{
Theorem \ref{thm:HB}  includes
the discovery that the property of ``uniform boundedness of multiplicity''
is determined
 only by the complexification $(G_{\mathbb{C}}, X_{\mathbb{C}})$
 and is independent of a real form $(G,X)$.
It is expected that this kind of statements could be generalized
 for reductive algebraic groups over non-archimedean local fields.
Recently, Sakellaridis--Venkatesh \cite{SaVe17} has obtained 
 some affirmative results in this direction.
}}
\end{remark}

\begin{remark}
\label{rem:Wh}
{\rm{
Theorems \ref{thm:HP} and \ref{thm:HB} give solutions 
 to Problem \ref{prob:6.2} (1) and (2), respectively.
More generally, these theorems hold not only for 
the space  $C^{\infty}(X)$ of functions
but also for the space of distributions and 
the space of sections of an equivariant vector bundle.
Furthermore, a generalization dropping the assumption 
  that the subgroup $H$ is reductive also holds, 
 see \cite[Thm.\ A, Thm.\ B]{xktoshima}
 for precise formulation.
For instance, the theory of the Whittaker model 
 considers the case
 where $H$ is a maximal unipotent subgroup, 
 see also Example \ref{ex:spherical} (5). 
Even for such a case
 a generalization of Theorems \ref{thm:HP} and \ref{thm:HB}
 can be applied.
}}
\end{remark}

\begin{remark}
{\rm{
We may also consider parabolic subgroups $Q$ 
 instead of a minimal parabolic subgroup.  
In this case, 
 we can also consider 
 \lq\lq{generalized Poisson transform}\rq\rq, 
 and extend the implication (i) $\Rightarrow$ (ii)
 in Theorem \ref{thm:HP}, 
 see \cite[Cor.~6.8]{xkProg2014}
 for a precise formulation.  
On the other hand, 
 an opposite implication (ii) $\Rightarrow$ (i)
 for parabolic subgroups $Q$ is not always true, 
see Tauchi \cite{Tauchi}.  
}}
\end{remark}

Theorems \ref{thm:HP} and \ref{thm:HB} suggest nice settings 
 of global analysis
 in which the \lq\lq{grip strength}\rq\rq\
 of representation theory is  \lq\lq{strong}\rq\rq.
The Whittaker model and the analysis on semisimple symmetric spaces
 may be thought of 
 in this framework
 as was seen in (6) and (1), 
respectively, 
 of Example \ref{ex:spherical}.  
As yet another set of problems, 
 let us discuss briefly the restriction 
 of representations to subgroups
 ({\bf{branching problems}}).

In the spirit of  \lq\lq{grip strength}\rq\rq\
 (Basic Problem \ref{prob:1.1}), 
 we may ask \lq\lq{grip strength
 of a subgroup}\rq\rq\ on an irreducible representation
 of a larger group
 as follows:

\begin{basicproblem}
[Grip strength in branching problem]
\label{pb:BP}
Let $\Pi$ be an irreducible representation of a group $G$.  
We regard $\Pi$ as a representation of a subgroup $G'$ by restriction, 
 and consider how many times another irreducible representation $\pi$ of $G'$
 occurs in the restriction $\Pi|_{G'}$:
\begin{enumerate}
\item[(1)] When is the multiplicity of every irreducible representation $\pi$ of $G'$ occurring in the restriction $\Pi|_{G'}$ finite?  
\item[(2)] When is the multiplicity of irreducible representation $\pi$ of $G'$ occurring in the restriction $\Pi|_{G'}$ uniformly bounded?  
\end{enumerate}
\end{basicproblem}

To be precise, 
 we need to clarify 
 what  \lq\lq{occur}\rq\rq\ means, 
 {\it{e.g., }}
 as a submodule, 
 as a quotient, 
 or as a support of the direct integral \eqref{eqn:MT}
 of the unitary representation, 
 {\it{etc}}.  
Furthermore, 
 since our concern is with infinite-dimensional irreducible representations, 
 the definition of \lq\lq{multiplicity}\rq\rq\
 depends also on the topology
 of the representation spaces of $\Pi$ of $G$ and $\pi$ of $G'$.  
Typical definitions of multiplicities include:
\begin{align}
\label{eqn:m1}
&\dim \invHom{G'}{\Pi^{\infty}|_{G'}}{\pi^{\infty}}, 
\\
\label{eqn:m2}
&\dim \invHom{G'}{\Pi|_{G'}}{\pi}, 
\\
\label{eqn:m3}
&\dim \invHom{G'}{\pi}{\Pi|_{G'}}, 
\\
\label{eqn:m4}
&\dim \invHom{{\mathfrak{g}}',K'}{\pi_{K'}}{\Pi_K}.   
\end{align}
Here $\Pi^{\infty}$, $\pi^{\infty}$ stand for smooth representations, 
 whereas $\pi_{K'}$ and $\Pi_K$ stand for the underlying $({\mathfrak{g}}',K')$-modules 
and $({\mathfrak{g}},K)$-modules.  
If $\Pi$ and $\pi$ are both unitary representations, 
 then the quantities \eqref{eqn:m2} and \eqref{eqn:m3} coincide.  
If \eqref{eqn:m4} $\ne 0$ in addition, 
 then all the quantities
 \eqref{eqn:m1}--\eqref{eqn:m4} coincide.  
In general the multiplicity 
 \eqref{eqn:m4} often vanishes, 
 and its criterion is given in \cite{TKAnn98, TKInvent98}.

Concerning the multiplicity \eqref{eqn:m1}, 
 see \cite{xkProg2014}, \cite{SunZhu12}
 and references therein for the general theory, 
 in particular, 
 for a geometric necessary and sufficient condition
 on the pair $(G, G')$
 such that (1) (or more strongly (2))
 of Basic Problem \ref{pb:BP} is always fulfilled.  
When the triple $(\Pi, G, G')$ satisfies
 finiteness
 (or more strongly, uniform boundedness)
 of the multiplicity
 in Basic Problem \ref{pb:BP}, 
 we could expect a detailed study
 of the restriction $\Pi|_{G'}$, 
 see \cite{xKVogan2015}, 
 for further \lq\lq{programs}\rq\rq\
 of branching problems of reductive groups, 
 such as the construction of \lq\lq{symmetry breaking operators}\rq\rq\
 and \lq\lq{holographic operators}\rq\rq\
 in concrete settings
\cite{KKP16, xtsbon, xksbonvec}.

%%%%%%%%%%%%%%%%%%%%%%%%%%%%%%%%%%%%%%%%%%
\section{Spectrum of the regular representation $L^2(X)$ ---a geometric criterion for temperdness}
\label{sec:tempered}
%%%%%%%%%%%%%%%%%%%%%%%%%%%%%%%%%%%%%%%%%

In the previous section,
we focused on \lq\lq{multiplicity}\rq\rq\ from the perspective
 of \lq\lq{grip strength}\rq\rq\ 
of a group on a function space
 and proposed (real) spherical varieties as
 \lq\lq{a nice framework for {\bf{detailed}} study of global analysis}\rq\rq.  
On the other hand, even in a case in which the \lq\lq{grip strength}\rq\rq\ 
of representation theory is \lq\lq{weak}\rq\rq,
 we may still expect to analyze the space of functions on $X$ from representation theory in a 
\lq\lq{{\bf{coarse}} standpoint}\rq\rq.  
In this section, 
 including {\it{non-spherical cases}},
let us focus on the support of the Plancherel measure
and consider the following problem.

Suppose that a Lie group $G$ acts on a manifold $X$
 with a Radon measure $\mu$
 and that $G$ leaves the measure invariant 
 so that $G$ acts naturally on the Hilbert space $L^2(X) \equiv L^2(X, d\mu)$
 as a unitary representation.  

\begin{basicproblem}
[Tempered space {\cite{BK2015}}]
\label{prob:6.7}
{\rm{
Find a necessary and sufficient condition
on a pair $(G,X)$ for which the regular representation $L^2(X)$ of $G$
is a tempered representation.
}}
\end{basicproblem}

We recall the general definition of tempered representations.  

\begin{defin}
[Tempered representation]
A unitary representation $\pi$ of a locally compact group $G$
 is called {\it{tempered}}
 if $\pi$ is weakly contained in $L^2(G)$, 
 namely,
 if any matrix coefficient $G \ni g \mapsto (\pi(g)u, v) \in {\mathbb{C}}$
 can be approximated
 by a sequence
 of linear combinations
 of matrix coefficients
 of the regular representation $L^2(G)$
 on every compact set of $G$.  
\end{defin}

The classification of {\it{irreducible}} tempered representations
 of real reductive linear Lie groups $G$ 
 was accomplished by Knapp--Zuckerman \cite{KZ}.  
In contrast to the long-standing problem
 of the classification of the unitary dual $\widehat G$, 
 irreducible tempered representations
 form a subset of $\widehat G$
 that is fairly well-understood.  
Loosely speaking, from the orbit philosophy
 due to Kirillov--Kostant--Duflo, 
 irreducible tempered representations are supposed to be
 obtained as a \lq\lq{geometric quantization}\rq\rq\
 of regular semisimple coadjoint orbits, 
 see {\it{e.g.,}} \cite{Krl2004, K94b}.

Tempered representations are unitary representations
 by definition, 
 however, 
 the classification theory of Knapp--Zuckerman played
 also a crucial role 
 in the Langlands classification
 of irreducible admissible representations
 (without asking
 if they are unitarizable or not)
 of real reductive Lie groups.

\vskip 1pc
The general theory of Mautner--Teleman tells
 that any unitary representation $\Pi$
 of a locally compact group $G$
 can be decomposed 
 into the direct integral
 of irreducible unitary representations:
\begin{equation}
\label{eqn:MT}
  \Pi \simeq \int^{\oplus} \pi_{\lambda} d \mu (\lambda).  
\end{equation}

Then the following equivalence (i) $\Leftrightarrow$ (ii) holds 
(\cite[Rem.~2.6]{BK2015}):
\begin{enumerate}
\item[(i)]
$\Pi$ is tempered;
\item[(ii)]
irreducible representation 
$\pi_{\lambda}$ is tempered
 for a.e. $\lambda$ with respect to the measure $\mu$.
\end{enumerate}

The irreducible decomposition of the regular representation
 of $G$ on $L^2(X)$ is called
 the Plancherel-type theorem for $X$.  
Thus, 
 if the Plancherel formula is  \lq\lq{known}\rq\rq, 
 then we should be able to answer
 Basic Problem \ref{prob:6.7}
 in principle. 
However, 
 things are not so easy:

\begin{observation}
\label{rem:6.15}
{\rm{The Plancherel-type theorem 
for semisimple symmetric spaces $G/H$
 was proved 
 by T.~Oshima, 
 P.~Delorme, 
 E.~van den Ban, 
 and H.~Schlichtkrull 
 (up to nonvanishing condition of discrete series representation 
 with singular parameters).  
However, 
 it seems that 
 a necessary and sufficient condition
on a symmetric pair $(G,H)$ for which $L^2(G/H)$ is tempered
had not been found until the general theory 
\cite{BK2015} 
 is established by a completely different approach.
In fact, 
 it is possible to show
 that temperedness of $L^2(G/H)$ implies
 a simple geometric condition 
 that $(G/H)_{\operatorname{Am}}$ is dense in $G/H$
 (see the second statement of Theorem \ref{thm:2.8}) from their Plancherel-type formula 
 in the case where $G/H$ is a symmetric space, 
 whereas there is a counterexample
 to the converse statement, 
 as was found 
 in \cite{BK2015}.  
If one employs the Plancherel-type formula
 in order to derive the right answer to Problem \ref{prob:6.7}
 for symmetric spaces $G/H$, 
 one will need a precise (non-)vanishing condition
on certain cohomologies 
 (Zuckerman derived functor modules)
 with singular parameters, 
 and such a condition 
is combinatorially complicated
 in many cases
(\cite{TK1992, Trapa2001}).  
}}
\end{observation}

\begin{observation}
{\rm{
More generally, when $X_{\mathbb{C}}$ 
is not necessarily a spherical variety of $G_{\mathbb{C}}$,
as shown in Theorem \ref{thm:HB},
the ring ${\mathbb{D}}_G(X)$ of $G$-invariant differential operators on $X$
is not commutative and so we cannot use effectively the existing method 
on non-commutative harmonic analysis
 based on an expansion of functions on $X$ into joint eigenfunctions
 with respect to the commutative ring ${\mathbb{D}}_G(X)$, 
 {\it{cf.}} Example \ref{ex:PDE}.
}}
\end{observation}

As observed above,
to tackle Basic Problem \ref{prob:6.7},
 one needs to develop a new method itself.
As a new approach,
 Benoist and I utilised an idea of dynamical system
 rather than differential equations.
We begin with some basic notion:

\begin{defin}
[Proper action]
\label{def:proper}
Suppose that a locally compact group $G$ acts 
continuously on a locally compact space $X$.  
This action is called {\bf{proper}}
 if the map 
\[
  G \times X \to X \times X, 
\qquad
 (g,x) \mapsto (x, g \cdot x)
\]
is proper, 
namely,
 if 
\[
  G_S:=\{s \in G: g S \cap S \ne \emptyset \}
\]
is compact for any compact subset $S$ of $X$.  
\end{defin}

If $G$ acts properly on $X$, 
 then the stabilizer of any point $x \in X$ in $G$ is compact.  
On the other hand, 
 if $H$ is a compact subgroup of $G$, then $L^2(G/H) \subset L^2(G)$ holds, 
 hence the regular representation on $L^2(G/H)$ is tempered.  
The following can be readily drawn from this.

\begin{example}
\label{ex:propertemp}
If the action of a group $G$ on $X$ is proper (Definition \ref{def:proper}),
then the regular representation in $L^2(X)$ is tempered.
\end{example}

Therefore, in the study of Basic Problem \ref{prob:6.7}, 
we focus on the nontrivial case that the action of $G$ on $X$ is not proper.
Properness of the action is {\bf{qualitative property}}, 
 namely, 
 there exists a compact subset $S$ of $X$ such that 
the set 
$
G_S=\{g \in G: g S \cap S \ne \emptyset\}
$
is noncompact.  
In order to shift it  {\bf{quantitatively}},
we consider the volume $\operatorname{vol}(g S \cap S)$.  
Viewed as a function on $G$, 
\begin{equation}
\label{eqn:2.1}
  G \ni g \mapsto \operatorname{vol}(g S \cap S) \in {\mathbb{R}}
\end{equation}
is a continuous function of $g \in G$.  
Definition \ref{def:proper} tells that the $G$-action on $X$
 is not proper
 if and only if 
 the support of the function \eqref{eqn:2.1}
 is noncompact for some compact subset $S$ of $X$.  
Hence the  \lq\lq{decay}\rq\rq\
 of the function \eqref{eqn:2.1} 
 at infinity may be considered
 as capturing quantitatively 
 a \lq\lq{degree}\rq\rq\ of non-properness of the action.
By pursuing this idea, Basic Problem \ref{prob:6.7} is settled in 
 Benoist--Kobayashi \cite{BK2015, BK2017}
 when $X$ is an algebraic $G$-variety for a reductive group $G$.  
To describe the solution,
 let us introduce a piecewise linear function 
 associated to a finite-dimensional representation of a Lie algebra.

\begin{defin}
For a representation 
$
   \sigma \colon {\mathfrak{h}} \to {\operatorname{End}}_{\mathbb{R}}(V)
$
of  a Lie algebra ${\mathfrak{h}}$ on a finite-dimensional real vector space $V$,
we define a function $\rho_V$
 on ${\mathfrak{h}}$ by
\begin{align*}
  \rho_V \colon {\mathfrak{h}} \to {\mathbb{R}},&
\quad
 \text{$Y \mapsto$
 the sum of the absolute values of the real parts}\\
 &\text{\hspace{0.4in} of the eigenvalues of $\sigma(Y)$ on
 $V \otimes_{\mathbb{R}} {\mathbb{C}}$.}
\end{align*}
\end{defin}

The function $\rho_V$ is uniquely determined by the restriction to 
a maximal abelian split subalgebra  ${\mathfrak{a}}$ of $\mathfrak{h}$.
Further, the restriction $\rho_V|_{\mathfrak{a}}$ is a piecewise linear function on ${\mathfrak{a}}$,
namely, 
 there exist finitely many convex polyhedral cones which cover
$\mathfrak{a}$ and on which $\rho_V$ is linear.

\begin{example}
When $(\sigma,V)$ is the adjoint representation $({\operatorname{ad}}, {\mathfrak{h}})$,
the restriction $\rho_{\mathfrak{h}}|_{\mathfrak{a}}$ can be computed by using
a root system.
It coincides with twice the usual \lq\lq{$\rho$}\rq\rq\
 in the dominant Weyl chamber.  
\end{example}

With this notation, 
one can describe a necessary and sufficient condition for Basic Problem \ref{prob:6.7}.

\begin{theorem}
[Criterion for temperedness of $L^2(X)$, {\cite{BK2017}}]
\label{thm:BK}
Let 
$G$ be a reductive Lie group and
$H$ a connected closed subgroup of $G$.
We denote by ${\mathfrak{g}}$ and ${\mathfrak{h}}$
 the Lie algebras of $G$ and $H$, respectively.
Then the following two conditions on a pair $(G,H)$ are equivalent.
\par\indent
{\rm{(i)}}
{\rm{(global analysis)}}\enspace
The regular representation $L^2(G/H)$ is tempered.
\par\indent
{\rm{(ii)}}
{\rm{(combinatorial geometry)}}\enspace
$\rho_{{\mathfrak{h}}} \le \rho_{{\mathfrak{g}}/{\mathfrak{h}}}$.  
\end{theorem}

\begin{remark}
\label{rem:6.14}
If 
$G$ is an algebraic group
 acting on an algebraic variety $X$,
then, even when $X$ is not a homogeneous space of $G$,
one can give an answer to Basic Problem \ref{prob:6.7} by 
applying Theorem \ref{thm:BK} to generic $G$-orbits
(\cite{BK2017}).
\end{remark}

Theorem \ref{thm:BK} was proved
 in \cite{BK2015}
 in the special case where both $G$ and $H$ are real algebraic 
 reductive groups.  
In this case, 
the following theorem also holds:

\begin{theorem}
[{\cite{BK2015}}]
\label{thm:2.8}
Let $G \supset H$ be a pair of real algebraic reductive Lie group.  
We set 
\begin{align*}
(G/H)_{\operatorname{Am}}
:=
&\{x \in G/H: \text{the stabilizer of $x$ in $H$ is  amenable}\}
\\
(G/H)_{\operatorname{Ab}}
:=
&\{x \in G/H: \text{the stabilizer of $x$ in $H$ is abelian}\}.  
\end{align*}
Then the following implications hold.
\begin{alignat*}{2}
&\text{geometry}\qquad
&&\text{$(G/H)_{\operatorname{Ab}}$ is dense in $G/H$}
\\
&
&&\hphantom{MMM}
\Downarrow
\\
&\text{representation}\qquad
&&\text{$L^2(G/H)$ is tempered}
\\
&
&&\hphantom{MMM}
\Downarrow
\\
&\text{geometry}
&&\text{$(G/H)_{\operatorname{Am}}$ is dense in $G/H$.  }
\end{alignat*}
\end{theorem}

Since a complex Lie group is amenable
 if and only if it is abelian, 
 Theorem \ref{thm:2.8} implies the following:

\begin{corollary}
The following conditions on a pair of complex reductive Lie groups $(G,H)$
 are equivalent:
\begin{enumerate}
\item[(i)]
$L^2(G/H)$ is tempered. 
\item[(ii)]
$(G/H)_{\operatorname{Ab}}$ is dense in $G/H$.  
\end{enumerate}
\end{corollary}

\medskip

{\bf Acknowledgements}:\enspace 
The author was partially supported by Grant-in-Aid for Scientific Research (A) (18H03669), Japan Society for the Promotion of Science.  
He also would like to thank Professor Vladimir Dobrev
 for his warm hospitality.

\end{document}